\documentclass{amsart}[11pt]

\usepackage{amscd,amssymb,latexsym,url,verbatim,graphicx,color}
\usepackage{hyperref}
\hypersetup{colorlinks=true,linkcolor=[rgb]{0,0,1},filecolor=[rgb]{0,0,1}}
\usepackage[capitalize,nameinlink,noabbrev]{cleveref}
\usepackage{tikz,tikz-cd}

\usetikzlibrary{decorations.pathmorphing}

\usepackage{cases,amsmath}

\usepackage{pxfonts}
\usepackage[euler-digits,small]{eulervm}

\usepackage{tikz,tikz-cd}

\usepackage[dvips]{epsfig}

\title{Configuration spaces of labeled points on a circle with two anchors}

\author{Dmitry N. Kozlov}

\address{Department of Mathematics, University of Bremen, 28334
  Bremen, Federal Republic of Germany.}
	
\email{dfk@math.uni-bremen.de}

\address{Okinawa Institute of Science and Technology Graduate University,
1919-1 Tancha, Onna-son, Kunigami-gun,
Okinawa, Japan.}


\newtheorem{theorem}{Theorem}[section]
\newtheorem{df}[theorem]{Definition}
\newtheorem{thm}[theorem]{Theorem} 

\newtheorem{prop}[theorem]{Proposition}
\newtheorem{lm}[theorem]{Lemma} 
\newtheorem{crl}[theorem]{Corollary}

\newcommand{\nin}{\noindent} 
\newcommand{\pr}{\nin{\bf Proof.} }

\newcommand{\bo}{\partial}

\newcommand{\clb}{{\mathcal B}}
\newcommand{\clc}{{\mathcal C}}

\newcommand{\clf}{{\mathcal F}}

\newcommand{\clk}{{\mathcal K}}

\newcommand{\clm}{{\mathcal M}}

\newcommand{\mcr}{{\mathcal R}}
\newcommand{\clr}{\mcr}
\newcommand{\cs}{{\mathcal S}}

\newcommand{\zz}{{\mathbb Z}}
\newcommand{\dz}{{\mathbb Z_2}}

\newcommand{\es}{\emptyset}

\newcommand{\ra}{\rightarrow}
\newcommand{\sm}{\setminus}
\newcommand{\supp}{\text{\rm supp}\,}

\newcommand{\im}{\textrm{Im}}

\newcommand{\mycap}[1]{\caption{#1}}

\newcommand{\myfig}[2]{
\begin{figure}[hbt]
\centering
\input{#1.pdf_t}
\mycap{#2.}
\label{fig:#1}
\end{figure}}

\newcommand{\myset}[2]{\{#1\,|\,#2\}}

\graphicspath{{./figs/}}

\numberwithin{equation}{section}
\numberwithin{figure}{section}
\numberwithin{table}{section}

\newcommand{\acs}[3]{\Sigma(#1,#2,#3)}

\newcommand{\cosp}[2]{\Sigma_{#1}(S^1,#2)}
\newcommand{\os}[1]{\cosp 2{#1}}

\newcommand {\x}[1]{\alpha(#1)}
\newcommand {\y}[1]{\beta(#1)}

\newcommand{\muminus}{\mu_-}
\newcommand{\mup}{\clm^{\uparrow}}
\newcommand{\mdown}{\clm^\downarrow}
\newcommand{\muplus}{\mu_+}
\newcommand{\crit}{Crit}
\newcommand{\mycut}{}
\newcommand{\thom}{\text{\tt{Hom}}}

\begin{document}

\begin{abstract}
In this paper we calculate the homology of configuration spaces of $n$
points on a circle, subject to the condition that two pre-determined
points are included in the configuration.

 We make use of discrete Morse theory both to determine the Betti
 numbers, as well as to provide an explicit combinatorial description
 of the bases both for homology and cohomology.
\end{abstract}

\maketitle

\section{Anchored configuration spaces}

\subsection{Definition}

Configuration spaces constitute a well-studied class of topological
spaces. Given an arbitrary topological space $X$ and a positive
integer $n$, a configuration space of $n$ labeled points in $X$
consists of all tuples $(x_1,\dots,x_n)$ of distinct points in $X$,
with the topology inherited from that of the direct product $X^n$.

In this paper, motivated by some problems in logistics, we study a
variation of this class.  The two main differences are that first, the
points $x_i$ do not need to be distinct, and second, that one
requires, that a certain pre-determined discrete set of points is
covered by each such tuple. The formal definition is as follows.

\begin{df} 
\label{df:acs}
Let $X$ be a non-empty topological space, let $S$ be a set of $m$
points in $X$, $m\geq 0$, and let $n$ be an arbitrary positive
integer.  An {\bf anchored configuration space}, denoted $\acs XSn$,
is defined as the subspace of the direct product $X^n$, consisting of
all tuples $(x_1,\dots,x_n)$, such that $S\subseteq\{x_1,\dots,x_n\}$.
\end{df}

For simplicity, in this paper we will assume that $X$ can be given a
structure of CW complex.  We make the following observations.

\begin{itemize}
\item When $m=0$, we simply have $\acs X{\es}n=X^n$. Therefore,
we can henceforth assume that $m\geq 1$.
\item If   $n<m$, we have $\acs XSn=\emptyset$
\item If $n=m$, the space $\acs XSn$ is a collection of $n!$ points
  equipped with discrete topology.
\end{itemize}

Furthermore, it is convenient to extend the \cref{df:acs}, and to declare
$\acs{X}{\emptyset}{0}$ to be a topological space consisting of a
single point.

\begin{df}
Let $n$ be a positive integer, and let $I$ be an arbitrary non-empty
set.  An {\bf $I$-partition} of $n$ is collection of nonnegative integers
$\{n_i\}_{i\in I}$, such that $\sum_{i\in I}n_i=n$.
\end{df}
Note, that all but finitely many of $n_i$ must be equal to $0$.  We
may think about an $I$-partition as a function $\varphi$ from $I$ to
the set of nonnegative integers, by setting $\varphi(i):=n_i$.

\subsection{The case  of the disconnected spaces}

It turns out, that one can always assume that $X$ is connected, due to
the following simple observation.

\begin{prop} 
\label{prop:disc}
Let $X$ be a topological space, let $S$ be a discrete subset of $X$
and let $n$ be a positive integer. Assume $X=\cup_{i\in I} X_i$ is a
decomposition into connected components.\footnote{Since $X$ is assumed
to be a CW complex, these are the same as its path-connected
components.} For all $i\in I$, set $S_i:=X_i\cap S$.

The space $\acs XSn$ is homeomorphic to the space 
\[\coprod_{\varphi}\prod_{i\in I}\acs{X_i}{S_i}{\varphi(i)},\]
where the disjoint union is taken over all $I$-partitions $\varphi$.
\end{prop}
\pr For any tuple of points $(x_1,\dots,x_n)\in X^n$, and for any
$i\in I$, set $\varphi(i):=|X_i\cap\{x_1,\dots,x_n\}|$. This gives a
decomposition of $\acs XSn$ into a disjoint union with components
indexed by $I$-partitions.  When such an $I$-partition is fixed,
choosing points within each $X_i$ is independent of each other, so we
get a direct product of spaces. Finally, each such space, up to
homeomorphism, depends only on the value of $\varphi(i)$, not on the
indices of points which landed in $X_i$, so it is homeomorphic to
$\acs{X_i}{S_i}{\varphi(i)}$.  \qed

\subsection{The case $n=m+1$}

Now that we know that we can assume that $X$ is connected, we can see
how far we can get with the cases when the difference between $m$ and
$n$ is small.  As mentioned, for $m=n$ we get a finite discrete space,
so the first interesting case is when $n=m+1$. The following graph
plays the key role in the topology of $\acs XS{m+1}$.

\begin{df} \label{df:kmkn}
Let $m$ be a positive integer.
The graph $\thom(K_m,K_{m+1})$ is defined as follows:
\begin{itemize}
\item the vertices of $\thom(K_m,K_{m+1})$ are indexed by injections
$f:[m]\hookrightarrow[m+1]$;
\item the vertices indexed by $f$ and $g$ are connected by an edge if
  the values of the functions $f$ and $g$ differ for precisely one
  argument.
\end{itemize}
\end{df}

The graph $\thom(K_m,K_{m+1})$ is $m$-regular. It has $(m+1)!$
vertices and $\frac{m}{2}(m+1)!$ edges.

There are different ways to view the graph $\thom(K_m,K_{m+1})$. For
instance, our notation comes from the fact that it is a special case
of a more general construction of the so-called $\thom$-complexes.
Given any two graphs $T$ and $G$, one can construction a
prodsimplicial complex $\thom(T,G)$, whose cells are
multihomomorphisms from $T$ to $G$. When the graphs $T$ and $G$ are
complete graphs with $m$ and $m+1$ vertices, we obtain the graph from
\cref{df:kmkn}. We refer to \cite{BaK06,BaK07,Ko07,Ko08} for precise
definitions and further background for this area.
 
From another angle, $\thom(K_m,K_{m+1})$ can be viewed as a Cayley
graph of the symmetric group $\cs_{m+1}$ for the system of generators
consisting of $m$ transpositions $\myset{(i,m+1)}{1\leq i\leq m}$; the
reader is referred to \cite{MKS} for background on Cayley graphs.

\begin{prop}
Let $X$ be a non-empty connected topological space, and let $S$ be a set of $m$
points in $X$, $m\geq 1$.
The space $\acs XS{m+1}$ is homotopy equivalent to a wedge of $(m+1)!$
copies of $X$ with the graph $\thom(K_m,K_{m+1})$.

In particular, the homology groups of $\acs XS{m+1}$ with integer coefficients
are given by the formula
\[H_k(\acs XS{m+1})\approx
\begin{cases}
\sum_{(m+1)!}H_k(X),&\text{ for } k\neq 1;\\
\zz^{\zeta_m}\oplus\sum_{(m+1)!}H_1(X),&\text{ for } k=1,
\end{cases}\]
where $\zeta_m=\frac{1}{2}(m+1)!(m-2)+1$.
\end{prop}
\pr Set $T:=\acs XS{m+1}$, and let $(x_1,\dots,x_{m+1})$ denote an
arbitrary point of $T$.  We know that at most one of the points
$x_1,\dots,x_{m+1}$ does not belong to $S$.  Write
$S=\{s_1,\dots,s_m\}$

Assume $f:[m]\hookrightarrow[m+1]$ is an arbitrary injection.  Set
\[X_f:=\myset{(x_1,\dots,x_{m+1})}{s_i=x_{f(i)},\text{ for all }1\leq i\leq m}.\] 
In other words, $X_f$ is the locus of those configurations where the
point occupying $s_i$ is $x_{f(i)}$, and for $r=[m+1]\sm\im f$, the
point $x_r$ can be chosen arbitrarily.  Clearly, $T=\cup_f X_f$, where
the union ranges over all such injections.  Each $X_f$ is isomorphic
to $X$.

Assume $f$ and $g$ are two such injections. If their values differ for
two or more arguments, the corresponding subspaces $X_f$ and $X_g$ do
not intersect. Assume they differ for exactly one argument, say
$f(k)\neq g(k)$, and $f(i)=g(i)$, for $i\neq k$. In particular,
$[m+1]\sm\im f=g(k)$ and $[m+1]\sm\im g=f(k)$.  Then $X_f$ intersects
with $X_k$ in the single point $p_{f,g}=(x_1,\dots,x_{m+1})$ determined by
$x_{g(k)}=s_k$ and $x_{f(i)}=s_i$, for all $i\in[m]$.
 
Inspired by homotopy colimits, we can now deform the space $\acs
XS{m+1}$ as follows. First, replace each intersection point $p_{f,g}$
by an interval $I_{f,g}$ , with one end attached to $X_f$ at its copy
of $p_{f,g}$, and the other end attached to $X_g$, again at the
respective copy of $p_{f,g}$. This gives us a homotopy equivalent
space where the $(m+1)!$ copies of $X$ are connected by intervals. We
have assumed that $X$ is connected and that it has CW structure, so it
must be path-connected.

Next, choose now in each $X_f$ an arbitrary base point $b_f$. Let the
endpoints of each interval $I_{f,g}$ slide inside the spaces $X_f$ and
$X_g$ to the respective base points. Again, this produces a homotopy
equivalent topological space $Y$. It can be given the following
description: take a certain connected graph $\Gamma$ and attach a~copy
of $X$ to each of its vertices. Since this graph is connected, we then
see that $Y$ is homotopy equivalent to the wedge of this graph with
$(m+1)!$ copies of $X$. The graph can then be replaced by a wedge of
circles.

To complete the proof, we need to understand the structure of
$\Gamma$. The vertices of $\Gamma$ are indexed by the injections
$f:[m]\hookrightarrow[m+1]$. These are connected by an edge if and
only if the spaces $X_f$ and $X_g$. The condition for that gives
precisely the graph $\thom(K_m,K_{m+1})$.
\qed


Note that even when $X$ is connected and $m=1$, different choices of
$S$ may result in non-homeomorphic topological spaces. Take, for instance, 
$X$ to be the unit interval $[0,1]$, and set $n:=2$. When $S=\{0\}$,
the space $\acs XS2$ is homeomorphic to a closed interval, whereas when
$S=\{1/2\}$, the space $\acs XS2$ is homeomorphic to the plus sign.



\subsection{Anchored configuration spaces of graphs}

The case we are mostly interested in is when $X$ is a graph.  It can
be thought of in terms of applications to logistics as follows.
Imagine, that we have $n$ unique resources, and that they need to be
distributed among $m$ locations.  Imagine furthermore, that the
locations, to which the resources are distributed, are connected by
a~graph network, and that each resource can be shifted from its
location to a neighboring one. Simultaneous multiple shifts of
different resources are allowed, as long as at any point of this
shifting procedure in each node there remain some resources, which are
not being moved. The spaces $\acs XSn$ are combinatorial cubical
complexes which model this situation by introducing higher-dimensional
parameter spaces which encode the interplay of such shifts.

The situation when $X$ is a tree, which for our purposes here is the
same as finite $1$-dimensional contractible simplicial complex, was
studied in \cite{Ko21}, where the spaces $\acs XS2$ were called the
{\it Stirling complexes}. The following result has been proved there.

\begin{thm} (\cite[Theorem 2.5]{Ko21}).
\label{thm:treemain}
Assume $X$ is a tree. Let $S$ be a subset of the set of the vertices
of $X$, $|S|=m\geq 2$, and let $n$ be an integer, $n\geq m$.  
The anchored configuration space $\acs{X}S{n}$ is homotopy
equivalent to a wedge of $(n-m)$-dimensional spheres.

Let $f(m,n)$ denote the number of these spheres. Then $f(m,n)$ is
given by the following formula
\begin{multline} \label{eq:fmn}
f(m,n)=(m-1)^n-\binom{m}{1}(m-2)^n+\binom{m}{2}(m-3)^n+\dots \\
+(-1)^{m-1}\binom{m}{m-3}2^n+(-1)^m\binom{m}{m-2}.
\end{multline}
\end{thm}

Note, that \cref{thm:treemain} says in particular that in this case the homotopy
type of $\acs{X}S{n}$ only depends on the cardinality of~$S$, not on the 
specific choice of the anchor points.

Once the case of trees (and, due to \cref{prop:disc}, of forests) is 
understood, it is natural to let $X$ be a circle graph. This is the same as to 
consider $\acs XSn$, when $X=S^1$. Clearly, in this setting
the space $\acs XSn$ is uniquely determined by $n$ {\it up to homeomorphism}. 

When the set $S$ consists of a single point, the space $\acs XSn$ is homotopy
equivalent to a punctured torus of dimension $n$. Indeed, it can be seen as 
a cubical complex obtained from the canonical cubical decomposition of 
the $n$-torus\footnote{View $S^1$ as a CW complex with a single vertex, and 
consider the direct product of these CW complexes}
by deleting the top-dimensional cube. Accordingly, the Betti numbers are 
given by $\beta_d(\acs XSn)=\binom{n}{d}$, for $d=0,\dots,n-1$.

The main result of this paper is determining the Betti numbers of 
$\acs XSn$, for the case when $|S|=2$.
In what follows, we refer to \cite{Fu,GH,Hat,Mu} for the background 
in algebraic topology.


\section{Cubical structure on $\os n$}

From now on, we shall exclusively deal with the case of configurations
of $n$ points on a circle, subject to the condition that two fixed
points must occur among the points in the configuration.  For any
$n\geq 2$, we let $\Omega_n$ denote this space $\acs XSn$, where
$|S|=2$.  Our purpose is to determine the Betti numbers
of~$\Omega_n$. 

By recording which of the points in the configuration can be found in
the anchor set, we see that the topological space $\Omega_n$ has a
canonical cubical complex structure. In this structure, the cubes are
indexed by the $4$-tuples $(A,B,C,D)$, where
\begin{enumerate}
\item[(1)] each $(A,B,C,D)$ is an ordered set partition of the set $[n]=\{1,\dots,n\}$;
\item[(2)] the sets $A$ and $C$ are non-empty.
\end{enumerate}
The sets $A$ and $C$ give us the set of points in the anchor set, while
the sets $B$ and $D$ tell us on which side of the anchor points lie
the remaining configuration points; see \cref{fig:f1}.

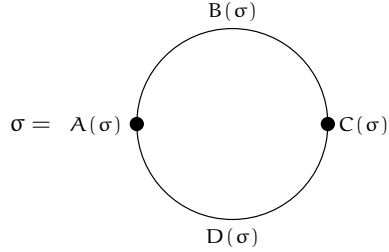
\myfig{f1}{Encoding of the cubes of $\Omega_n$}

When $\sigma$ is a cube of $\Omega_n$, and $\sigma=(A,B,C,D)$, to
express the dependence on $\sigma$, we shall also write
$\sigma=(A(\sigma),B(\sigma),C(\sigma),D(\sigma))$.  This way we shall
be able to refer to single parts of $\sigma$ when necessary.

For brevity we use the following notations:
\begin{itemize}
\item for $S\subseteq[n]$ and $x\in[n]\sm S$, we let $S+x$ denote $S\cup x$;
\item for $S\subseteq[n]$ and $x\in S$, we let $S-x$ denote $S\sm x$;
\item for $S\subseteq[n]$ and $x\in[n]\sm S$, we write $x<S$, if $x<y$,
for all $y\in S$; 
\newline same way we write $x>S$, if $x>y$, for all $y\in S$.
\end{itemize}
Note, that these abbreviations do not override the regular
set-theoretical notations, so, when convenient, we shall still use the
latter ones.

The cubical boundary operator over $\dz$ is given by
\begin{multline}
\label{eq:cbo}
\bo(A,B,C,D)=\sum_{x\in B}(A+x,B-x,C,D)+\sum_{x\in B}(A,B-x,C+x,D)\\
+\sum_{x\in D}(A+x,B,C,D-x)+\sum_{x\in D}(A,B,C+x,D-x).
\end{multline}

The dimension of the cube indexed by $(A,B,C,D)$ is equal to
$|B|+|D|$.  Accordingly, for $0\leq d\leq n-2$, the number of
$d$-cubes in $\Omega_n$ is equal to $\binom nd 2^d(2^{n-d}-2)=\binom
nd(2^n-2^{d+1})$, so its $f$-vector is
\[f(\Omega_n)=(2^n-2,n(2^n-4),\binom n2(2^n-8),\dots,\binom n3(2^n-2^{n-2}),\binom n2 2^{n-1}).\]

This gives us the following formula for the Euler characteristic of
$\Omega_n$: 
\begin{multline*}
\mu(\Omega_n)=\sum_{d=0}^{n-2}(-1)^d \binom nd(2^n-2^{d+1})=
2^n\sum_{d=0}^{n-2}(-1)^d \binom nd-2\sum_{d=0}^{n-2}(-1)^d \binom nd 2^d
\\=2^n(-1)^n(n-1)-2(-1)^n(2^{n-1}n+1-2^n)=(-1)^n(2^n-2).
\end{multline*}
This formula will be confirmed when we have calculated the Betti
numbers of~$\Omega_n$.

\mycut

\section{Discrete Morse theory for chain complexes}
\label{sect:admt}

Discrete Morse theory, \cite{Kn,Ko08,Ko20b,Sc19}, has emerged in many 
situations in combinatorial topology, \cite{Ca,EH}, 
as a useful tool for calculating homology groups, and perhaps
even understanding the homotopy type of the involved spaces. It has
been argued in \cite{admt} and \cite[Section 11.3]{Ko08} that an
algebraic point of view on this theory is rather beneficial. We now
sketch the basic tenets of that approach as applied to our context.

\begin{df}
Let $\clc=(C_*,\bo_*)$ be a chain complex of vector spaces over $\dz$, and
let $\{\Gamma_d\}_d$  be a collection of sets, such that
$\Gamma_d$ is a basis of $C_d$, for all $d$.
The union $\Gamma=\cup_d\Gamma_d$ is called a {\bf basis of $\clc$}.
\end{df}

Assume now we are given a chain complex $\clc$ over $\dz$ and a basis
$\Gamma$. We set $\Gamma_d:=\Gamma\cap C_d$, for all $d$.  Each chain
$\sigma\in C_d$ is a sum of some of the elements of $\Gamma_d$.  The
set of these elements is called the support of $\sigma$, and is
denoted by $\supp\sigma$, so we have
$\sigma=\sum_{\tau\in\supp\sigma}\tau$, for all~$C_d$.

\begin{df}
When $\sigma,\tau\in\Gamma$, we say that $\sigma$ {\bf covers} $\tau$
if $\tau$ is contained in the support of~$\sigma$. In such a case, we
write $\sigma\succ\tau$.
\end{df}

Obviously, when $\sigma\in\Gamma_d$ and $\sigma\succ\tau$, we must
have $\tau\in\Gamma_{d-1}$. Reversely, we shall say $\tau$ {\it is
  covered} by $\sigma$, and write $\tau\prec\sigma$.

\begin{df}
An {\bf acyclic matching} on $\Gamma$ is an
involution $\mu:\clm\ra\clm$, where $\clm$ is some subset of
$\Gamma$, such that
\begin{enumerate}
\item[(1)] for each $\alpha\in\clm$, the element $\mu(\alpha)$ either
  covers $\alpha$, or is covered by~$\alpha$;
\item[(2)] there do not exist distinct $a_1,\dots,a_t\in\Gamma$, 
such that $t\geq 2$, and 
\begin{equation}\label{eq:acycle}
\mu(a_t)\succ a_1\prec\mu(a_1)\succ a_2\prec\mu(a_2)\succ\dots
\prec \mu(a_{t-1})\succ a_t\prec \mu(a_t).
\end{equation}
\end{enumerate}
\end{df}

We decompose $\clm=\mup\cup\mdown$, where $\mup$ consists of all
$\alpha$, such that $\mu(\alpha)$ is covered by $\alpha$ in $\Gamma$,
while $\mdown$ consists of all $\alpha$, such that $\mu(\alpha)$
covers $\alpha$ in $\Gamma$. When we need to be specific, we also
write $\muplus(\alpha)$ instead of $\mu(\alpha)$, for
$\alpha\in\mdown$, and we write $\muminus(\alpha)$, when
$\mu\in\mup$. Furthermore, we let $\crit$ denote the complement
of~$\clm$, and set $\crit_d:=\Gamma_d\cap\crit$, for all~$d$. We say
that the basis elements in $\clm$ are \emph{matched} and the basis
elements in $\crit$ are \emph{critical}.

\begin{df}
Let $\sigma$ be a critical basis element. An {\bf alternating path} starting at
$\sigma$ is any sequence of basis elements
\begin{equation}
\label{eq:altpath}
\sigma\succ\tau_1\prec\muplus(\tau_1)\succ\dots\succ\tau_q\prec\muplus(\tau_q)\succ\tau,
\end{equation}
where  $q$ is a nonnegative integer and $\tau$ is a critical basis element.
\end{df}

Given an alternating path \eqref{eq:altpath}, we set $p_\bullet:=\sigma$ and 
$p^\bullet:=\tau$.

\begin{df}
\label{df:mc}
Assume $\clc$ is a chain complex over $\dz$ and assume $\Gamma$ is
a~basis of $\clc$, with an acyclic matching $\mu:\clm\ra\clm$. The
{\bf Morse complex} $\clc^\clm=(C_*^\clm,\bo_*^\clm)$ is defined as
follows:
\begin{itemize}
\item for every $d$, $C_d^\clm$ is a vector space over $\dz$, with a
  basis $\{c_\sigma^\clm\}_{\sigma\in\crit_d}$ indexed by the critical
  basis elements of dimension~$d$;
\item for each $\sigma\in\crit_d$ the value of the boundary operator
  on $c_\sigma^\clm$ is given by
\[\bo_d^\clm(c_\sigma^\clm):=\sum_p c_{p_\bullet}^\clm,\] 
where the sum is taken over all alternating paths $p$ satisfying
$p^\bullet=\sigma$.
\end{itemize}
\end{df}

It is a~known fact, see, e.g., \cite[Section 11.3]{Ko08},
\cite[Chapter 15]{Ko20b}, that $(C_*^\clm,\bo_*^\clm)$ is a
well-defined chain complex.


We now proceed with the core of the algebraic discrete Morse
theory. The following sequence of statements is derived from 
\cite[Theorem 11.24]{Ko08} and \cite[Chapter 15]{Ko20b}.

Assume, as above, that $\clc$ is a~chain complex over $\dz$, $\Gamma$
is a basis of $\clc$, and $\mu:\clm\ra\clm$ is an acyclic
matching. Each $C_d$ is then given a~basis, which can be seen as
a~disjoint union of the following three sets: $\crit_d$,
$\mup_d:=\mup\cap\Gamma_d$, and $\mdown_d:=\mdown\cap\Gamma_d$.  The
next theorem produces a~\emph{new basis} with the same indexing sets,
but with much improved boundary values, so that the chain complex
splits as a~direct sum.

\begin{thm} 
\label{thm:admt}
For each $d$,  there exist sets of vectors
\[B_d=\{b_\sigma\}_{\sigma\in\crit_d},\quad
U_d=\{u_\tau\}_{\tau\in\mup_d},\quad
L_d=\{l_\rho\}_{\rho\in\mdown_d},\]
such that the following statements hold:
\begin{enumerate}
\item[(1)] $B_d\cup U_d\cup L_d$ is a basis of $C_d$, for all~$d$;
\item[(2)] $\supp b_\sigma\cap \crit_d=\sigma$, for all $\sigma\in\crit_d$;
\item[(3)] $\bo b_\sigma=\sum_p b_{p_\bullet}$, where the sum is taken over all
alternating paths $p$ satisfying $p^\bullet=\sigma$;
\item[(4)] $\bo u_\tau=l_{\muminus(\tau)}$,  for all $\tau\in\mup_d$;
\item[(5)] $\bo l_\rho=0$, for all $\rho\in\mdown_d$.
\end{enumerate}
\end{thm}

\pr See \cite[Theorem 11.24]{Ko08} or \cite[Chapter 15]{Ko20b}. Note,
how the proof there follows a procedure closely resembling the
Gaussian elimination. \qed

\begin{crl} \label{crl:1}
The chain complex $\clc$ splits as a~direct sum $\clb\oplus\clr$ of
chain complexes, where $\clb$ is isomorphic to the Morse chain complex
$\clc^\clm$, and $\clr$ is acyclic.\footnote{Note the different use of
the word \emph{acyclic}: here it means that the homology groups of
$\clr$ vanish.}

In particular, we have $H_*(\clc)\cong H_*(\clc^\clm)$.
\end{crl}

\Cref{crl:1} says that the Morse complex calculates the homology of
the original complex $\clc$. When the Morse boundary operator
vanishes, much more detailed information can be obtained about the
homology and the cohomology generators.

Recall that the cochains can be seen as functions which can be
evaluated on chains. Furthermore, for a basis element
$\sigma\in\Gamma_d$, we let $\sigma^*$ denote the cochain which
evaluates to $1$ on $\sigma$, and evaluates to $0$ on any other
element of $\Gamma_a$. We call $\sigma^*$ the {\it dual} of~$\sigma$.

\begin{thm}\label{thm:admt2}
Assume that the boundary operator in the Morse complex $\clc^\clm$
vanishes.  Then, the following hold.
\begin{enumerate}
\item[(1)] The chains $\{b_\sigma\}_{\sigma\in\crit_d}$ from
  \cref{thm:admt} form a basis for $H_d(\clc)$.
\item[(2)] The cochains $\{\sigma^*\}_{\sigma\in\crit_d}$ form a basis for
  $H^d(\clc)$, for all~$d$.
\item[(3)] Any collection of cycles
  $\{d_{\sigma}\}_{\sigma\in\crit_d}$, such that $\sigma$ is the
  unique critical basis element in the support of $d_\sigma$, for all
  $\sigma\in\crit_d$, forms a basis for $H_d(\clc)$.
\end{enumerate} 
\end{thm}
\pr Statement (1) is immediate. 

To see (2), note that for every $\sigma\in\crit_d$, the cochain
$\sigma^*$ evaluates to $1$ on $b_\sigma$, and evaluates to $0$ on
$b_\tau$, for any $\tau\in\crit_d$, $\tau\neq\sigma$. Therefore,
$\{\sigma^*\}_{\sigma\in\crit_d}$ is a basis for $H^d(\clc)$.

Finally, to see (3), note that for every $\sigma\in\crit_d$, the
cochain $\sigma^*$ evaluates to $1$ on $d_\sigma$, while for any other
$\tau\in\crit_d$, $\tau\neq\sigma$, the cochain $\tau^*$ evaluates to
$0$ on~$d_\sigma$. \qed

The algebraic framework of chain complexes can readily be adapted to
cubical complexes.  So assume $\clk$ is a cubical complex. For $0\leq
d \leq \dim\clk$, we let $\clk(d)$ denote the set of all
$d$-dimensional cubes in $\clk$. We consider the cubical chain complex
of $\clk$ over $\dz$, $\clc:=C_*(\clk;\dz)$. In $\clc$, each chain
group $C_d(\clk;\dz)$ is a vector space over $\dz$ with basis
$\clk(d)$.  We can then apply the algebraic discrete Morse theory, by
taking $\Gamma_d:=\clk(d)$, for all $d$. Note that here the covering
relation comes from the following combinatorial framework.

\begin{df}
For a cubical complex $\clk$, we let $\clf(\clk)$ denote its
{\bf face poset}.  Its elements are the non-empty cubes of $\clk$,
which are partially ordered by inclusion.
\end{df}

The covering relation in the partial order on $\clf(\clk)$ is denoted
by~$\succ$. The acyclicity  condition can then be formulated using
this combinatorial notion of covering.

\mycut

\section{Acyclic matching}

\subsection{Defining a matching on $\clf(\Omega_n)$}

We proceed to describe a specific matching rule for the cubical
complex $\Omega_n$.  From now on, we set $\clk:=\Omega_n$.

First, we define two pivot functions
\[\alpha,\beta:\clf(\clk)\ra [n],\] 
as follows:
\begin{align*}
\x\sigma&:=\min(A(\sigma)+ B(\sigma))\\
\y\sigma&:=\max(B(\sigma)+ C(\sigma)) .
\end{align*}  Clearly,
either $\x\sigma\in A$ or $\x\sigma\in B$, and, similarly, either
$\y\sigma\in B$ or $\y\sigma\in C$.


Now, for all $0\leq d\leq n-2$, we define the following pairs of sets:
\begin{align*}
\mup_1(d)&:=\myset{\sigma\in\Omega_n(d)}{\x\sigma\in B(\sigma)}, \\
\mdown_1(d)&:=\myset{\sigma\in\Omega_n(d)}{\x\sigma\in A(\sigma)\text{ and }|A(\sigma)|\geq 2}.
\end{align*}
and
\begin{align*}
\mup_2(d)&:=\myset{\sigma\in\Omega_n(d)}{A(\sigma)=\x\sigma
\text{ and }\y\sigma\in B(\sigma)}, \\
\mdown_2(d)&:=\myset{\sigma\in\Omega_n(d)}{A(\sigma)=\x\sigma,
\y\sigma\in C(\sigma),\,|C(\sigma)|\geq 2\text{, and }\y\sigma>\x\sigma}.
\end{align*}
The reason for our choice of notations will become clear shortly. For
now, note that the four sets $\mup_1(d)$, $\mdown_1(d)$, $\mup_2(d)$,
and $\mdown_2(d)$ are disjoint for all~$d$.


For any $0\leq d\leq n-3$, we define functions 
\begin{align*}
f:&\mdown_1(d)\ra\mup_1(d+1),\\
g:&\mdown_2(d)\ra\mup_2(d+1),
\end{align*}
by setting
\begin{align*}
f(\sigma)&:=(A(\sigma)-\x\sigma,B(\sigma)+\x\sigma,C(\sigma),D(\sigma)), \text{ and} \\
g(\sigma)&:=(A(\sigma),B(\sigma)+\y\sigma,C(\sigma)-\y\sigma,D(\sigma)),
\end{align*}
so, in words, the function $f$ moves $\alpha(\sigma)$ from $A$ to $B$,
while the function $g$ moves $\beta(\sigma)$ from $C$ to $B$. Note,
that $f(\sigma)$ covers $\sigma$ and $g(\sigma)$ covers $\sigma$, for all 
$\sigma$'s in the definition domain of the respective function.

\begin{prop} \label{prop:fg}
The functions $f$ and $g$ are well-defined bijections, whose inverses are given by 
\begin{align*}
f^{-1}(\sigma)&:=(A(\sigma)+\x\sigma,B(\sigma)-\x\sigma,C(\sigma),D(\sigma)), \text{ and} \\
g^{-1}(\sigma)&:=(A(\sigma),B(\sigma)-\y\sigma,C(\sigma)+\y\sigma,D(\sigma)).
\end{align*}
\end{prop}
\pr Let us start with $f$. If $\sigma\in\mdown_1(d)$, then
$|A(\sigma)|\geq 2$, so $A(\sigma)-\x\sigma\neq\es$. This means that
$f(\sigma)\in\clk(d+1)$.  Furthermore, since $A(\sigma)+B(\sigma)
=A(f(\sigma))+B(f(\sigma))$, we have $\x\sigma=\x{f(\sigma)}$. This
means that $\x{f(\sigma)}\in B(f(\sigma))$, so
$f(\sigma)\in\mup_1(d+1)$, and hence $f$ is well-defined. Finally,
again since $\x\sigma=\x{f(\sigma)}$, the inverse of $f$ is as stated
in the formulation of the proposition.

Now, let us consider the function $g$.  Obviously, the functions $g$
and $g^{-1}$, as the latter is defined in the proposition, are
inverses of each other, as long as they are well-defined.

First, we see that $g$ is well-defined.  Take
$\sigma\in\mdown_2(d)$. We have $|C(\sigma)|\geq 2$, so
$C(\sigma)-\y\sigma\neq\es$, and therefore $g(\sigma)\in\clk(d+1)$.
We have $B(\sigma)+C(\sigma)=B(g(\sigma))+C(g(\sigma))$, so
$\y\sigma=\y{g(\sigma)}$, which implies $\y{g(\sigma)}\in
B(g(\sigma))$. 

Clearly, $A(g(\sigma))=A(\sigma)=\x\sigma$. Furthermore, 
\[\x{g(\sigma)}=\min(\x\sigma+B(\sigma)+\y\sigma)=\x\sigma,\]
since $\x\sigma<\y\sigma$ and $\x\sigma<B(\sigma)$.  It follows that
$g(\sigma)\in\mup_2(d+1)$.

Finally, let us see that $g^{-1}$ is well-defined. Take
$\sigma\in\mup_2(d)$.  It is obvious that
$g^{-1}(\sigma)\in\clk(d-1)$, and that $|C(g^{-1}(\sigma))|\geq 2$.
Again, since
$B(\sigma)+C(\sigma)=B(g^{-1}(\sigma))+C(g^{-1}(\sigma))$, we have
$\y\sigma=\y{g^{-1}(\sigma)}$, and so $\y{g^{-1}(\sigma)}\in
C(g^{-1}(\sigma))$.

Similar to the above,
$\x{g^{-1}(\sigma)}=\min(\x\sigma+B(\sigma)-\y\sigma)=\x\sigma$, since
$\x\sigma<\y\sigma$ and $\x\sigma<B(\sigma)$. It follows, that
$A(g^{-1}(\sigma))=\x{g^{-1}(\sigma)}$, and that
$\x{g^{-1}(\sigma)}<\y{g^{-1}(\sigma)}$. Thus
$g^{-1}(\sigma)\in\mdown_2(d)$, and the proof is finished. \qed

We now set 
\begin{align*}
\mdown(d)&:=\mdown_1(d)\cup\mdown_2(d),\text{ for all }d=0,\dots,n-3;\\
\mup(d)&:=\mup_1(d)\cup\mup_2(d),\text{ for all }d=1,\dots,n-2;
\end{align*}
\[\mup_1:=\cup_{d=1}^{n-2}\mup_1(d),\quad \mup_2:=\cup_{d=1}^{n-2}\mup_2(d);\]
\[\mdown_1:=\cup_{d=0}^{n-3}\mdown_1(d),\quad \mdown_2:=\cup_{d=0}^{n-3}\mdown_2(d);\]
\[\mdown:=\mdown_1\cup\mdown_2=\cup_{d=0}^{n-3}\mdown(d),\quad 
\mup:=\mup_1\cup\mup_2=\cup_{d=1}^{n-2}\mup(d);\]
\[\clm_1:=\mup_1\cup\mdown_1,\quad\clm_2:=\mup_2\cup\mdown_2;\]
\[\clm:=\mdown\cup\mup=\clm_1\cup\clm_2.\]
By \cref{prop:fg}, we know that functions $f$ and $g$ define a bijection from
$\mdown$ to $\mup$, whose inverse is given by the combination of the inverses
of $f$ and of~$g$. We denote this bijection by $\muplus$, its inverse by $\muminus$,
and the resulting involution of $\clm$ simply by~$\mu$.


\subsection{Critical cubes of $\clm$}

To describe the set of critical cubes of $\clm$, let us define the following three sets
\begin{align*}
C_1&:=\myset{\sigma\in\clf(\Omega_n)}{A(\sigma)=\x\sigma\text{, }
B(\sigma)=\emptyset\text{, }|C(\sigma)|\geq 2
\text{, and }\y\sigma<\x\sigma},\\
C_2&:=\myset{\sigma\in\clf(\Omega_n)}{A(\sigma)=\x\sigma\text{, }B(\sigma)=\emptyset
\text{, and }C(\sigma)=\y\sigma},\\
C_3&:=\myset{\sigma\in\clf(\Omega_n)}{A(\sigma)=\x\sigma\text{, }C(\sigma)=\y\sigma
\text{, } B(\sigma)\neq\emptyset
\text{, and }\x\sigma<B(\sigma)<\y\sigma},
\end{align*}
see \cref{fig:crit}.

Note that all the cubes in $C_2\cup C_3$ have dimension $n-2$, whereas
each cube $\sigma\in C_1$ has dimension $0\leq |D(\sigma)|\leq n-3$.

As in the general case, let $\crit$ denote the set of the critical cubes with
respect to the matching $\clm$.

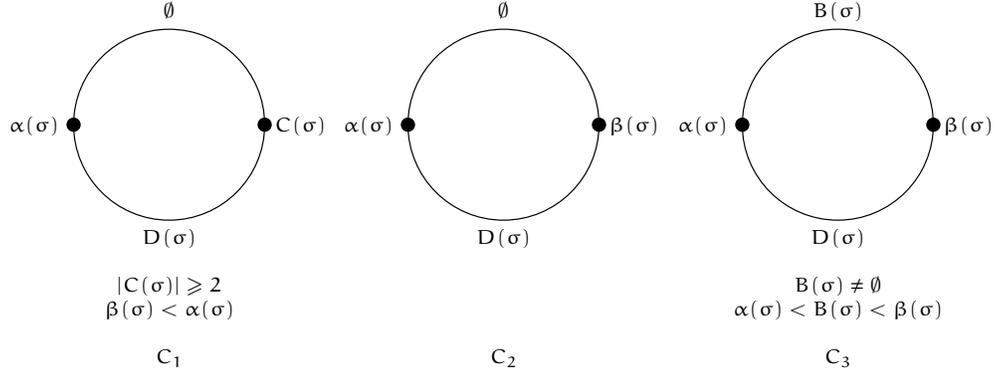
\myfig{crit}{The 3 types of critical cells with respect to $\clm$}

\begin{prop}
We have $\crit=C_1\cup C_2\cup C_3$.
\end{prop}
\pr Clearly, any critical cube $\sigma$ must satisfy
$A(\sigma)=\x\sigma$, or else it would belong to the set
$\clm_1$. This is satisfied for all $\sigma\in C_1\cup C_2\cup C_3$,
so we only need to consider $\sigma$'s, for which
$A(\sigma)=\x\sigma$.

Let us now look at those $\sigma$, for which
$B(\sigma)=\emptyset$. Such a cube is critical if and only if it does
not belong to $\mdown_2$. This is the case if and only if at least one
of the following conditions is satisfied:
\begin{enumerate} 
\item[(1)] either $C(\sigma)=\y\sigma$,
\item[(2)] or $C(\sigma)-\y\sigma\neq\emptyset$, but
  $\y\sigma<\x\sigma$.
\end{enumerate}
The first case describes the set $C_2$, whereas the second case
describes $C_1$.

Now look at those $\sigma$ for which $B(\sigma)\neq\emptyset$.  By
construction, we have $\x\sigma<B(\sigma)$, so $\x\sigma<\y\sigma$.
If $\y\sigma\in B(\sigma)$, we have $\sigma\in\mup_2$.  If
$\y\sigma\in C(\sigma)$, but $|C(\sigma)|\geq 2$, we have
$\sigma\in\mdown_2$.  So, the only option remaining is that
$\y\sigma=C(\sigma)$, in which case $\sigma$ is critical, and we have
precisely described the set~$C_3$.  \qed

\begin{prop}
\label{prop:count}
For $d=0,\dots,n-3$, the number of critical cubes of dimension $d$ is given by $\binom{n}{d}$.
The number of critical cubes of dimension $n-2$ is equal to $2^n+\binom{n-1}{2}-2$.
\end{prop}
\pr For $d=0,\dots,n-3$, the critical cubes of dimension $d$ are precisely those
cubes in $C_1$, for which $|D(\sigma)|=d$. Since $\x\sigma>C(\sigma)$, the choice of $D(\sigma)$
determines the cube uniquely. Therefore, the number of such cubes is simply $\binom{n}{d}$.

The number of critical cubes of dimension $n-2$ is equal to $|C_2|+|C_3|$.
For $\sigma\in C_2$, the numbers $\x\sigma$ and $\y\sigma$ can be chosen independently, so
$|C_2|=n(n-1)$. For $\sigma\in C_3$ again the choice of $D(\sigma)$ determines $\sigma$
uniquely, since in this case $\x\sigma=\min([n]\sm D(\sigma))$ and 
$\y\sigma=\max([n]\sm D(\sigma))$. The only condition is that $|D(\sigma)|\leq n-3$.
It follows that $|C_3|=2^n-1-n-\binom{n}{2}$. Adding this to $n(n-1)$ yields the desired formula. 
\qed


\subsection{The acyclicity of $\clm$}

Let us now show that our matching satisfies the key property required
by discrete Morse theory.

\begin{prop} \label{prop:clmacyc}
The matching $\clm$ is acyclic.
\end{prop}
\pr
Assume $\clm$ is not acyclic, and pick a cycle 
\begin{equation}
\label{eq:c1}
\sigma_1\prec\muplus(\sigma_1)\succ\sigma_2\prec\muplus(\sigma_2)\succ\dots
\succ\sigma_q\prec\muplus(\sigma_q)\succ\sigma_1,
\end{equation}
where $q\geq 2$, and $\sigma_1,\dots,\sigma_q$ are distinct cubes of
dimension $d$, where $0\leq d\leq n-3$.  We traverse this cycle from
left to right. For convenience, set $\sigma_{q+1}:=\sigma_1$.

By definition of the matching, we have
$D(\sigma_i)=D(\muplus(\sigma_i))$, for $1\leq i\leq q$.  On the other
hand, since $\muplus(\sigma_i)\succ\sigma_{i+1}$, we have
$D(\muplus(\sigma_i))\subseteq D(\sigma_{i+1})$, for all $i$.
Therefore, we conclude
\[D(\sigma_1)\subseteq D(\sigma_2)\subseteq\dots
\subseteq D(\sigma_q)\subseteq D(\sigma_1),\]
which of course implies
\[D(\sigma_1)=\dots=D(\sigma_q)=D(\muplus(\sigma_1))=\dots
=D(\muplus(\sigma_q)).\]

Let us now see what happens to the structure of the cube labels
when one follows the cycle.

Consider an arbitrary index $i$, such that $\sigma_i\in\mdown_1$.  For
brevity, say $\sigma_i=(A,B,C,D)$.  We have
\[\muplus(\sigma_i)=(A-\x{\sigma_i},B+\x{\sigma_i},C,D).\] 
What are the options for $\sigma_{i+1}$? Assume 
\[\sigma_{i+1}=(A-\x{\sigma_i}+y,B+\x{\sigma_i}-y,C,D).\]
Since $\sigma_i\neq\sigma_{i+1}$, we have $y\neq\x{\sigma_i}$. Then
$\x{\sigma_{i+1}}=\x{\sigma_i}$, and, since $\x{\sigma_i}\in
B(\sigma_{i+1})$, we have $\sigma_{i+1}\in\mup_1$, which is a
contradiction.  We must therefore have
\[\sigma_{i+1}=(A-\x{\sigma_i},B+\x{\sigma_i}-y,C+y,D),\]
for some $y\in B+\x{\sigma_i}$. If $y\neq\x{\sigma_i}$, then again
$\x{\sigma_{i+1}}=\x{\sigma_i}$, this time since
\[\x{\sigma_i}\in A(\sigma_{i+1})+B(\sigma_{i+1})\subseteq A(\sigma_i)+B(\sigma_i).\]
This means that $\sigma_{i+1}\in\mup_1$, again leading to a
contradiction.  We therefore conclude that the only option is
\[\sigma_{i+1}=(A-\x{\sigma_i},B,C+\x{\sigma_i},D).\]

Let us now consider the case $\sigma_i\in\mdown_2$, so we can write
$\sigma_i=(\x{\sigma_i},B,C,D)$, with $|C|\geq 2$, $\y{\sigma_i}\in
C$, and $\y{\sigma_i}>\x{\sigma_i}$. By construction
$\x{\sigma_i}<B$. We then have
\[\muplus(\sigma_i)=(\x{\sigma_i},B+\y{\sigma_i},C-\y{\sigma_i},D).\]
Assume first
\[\sigma_{i+1}=(\x{\sigma_i},B+\y{\sigma_i}-x,C-\y{\sigma_i}+x,D).\]
Since $\sigma_i\neq\sigma_{i+1}$ we have $x\neq\y{\sigma_i}$.  We then
have $\x{\sigma_i}<B+\y{\sigma_i}-x$, because
$\x{\sigma_i}<\y{\sigma_i}$ and $\x{\sigma_i}<B$. This means that
$\x{\sigma_{i+1}}=\x{\sigma_i}$.  On the other hand,
$\y{\sigma_i}=\y{\sigma_{i+1}}$, so $\sigma_{i+1}\in\mup_2$, yielding
a contradiction. We can therefore conclude that in this case
\[\sigma_{i+1}=(A+x,B-x+\x{\sigma_i},C-\x{\sigma_i},D),\]
where $x$ may or may not be equal to $\x{\sigma_i}$.

These considerations imply that the matchings of the two types
alternate along the cycle \eqref{eq:c1}. We can therefore assume that
$q=2t$ for some $t\geq 1$, and that
$\sigma_1,\sigma_3,\dots,\sigma_{2t-1}\in\mdown_1$, whereas
$\sigma_2,\sigma_4,\dots,\sigma_{2t}\in\mdown_2$.


To finish the argument let us look closely at each of these
alternating steps.  Pick $1\leq k\leq t$, and consider
$\sigma_{2k-1}$, say $\sigma_{2k-1}=(A,B,C,D)$.  The argument that
follows is illustrated by \cref{fig:acyc} for the case $k=1$.

Since $\sigma_{2k-1}\in\mdown_1$, we have $|A|\geq 2$ and
$\x{\sigma_{2k-1}}\in A$. Accordingly, by our argument above, the
cycle continues with
\begin{align*}
\muplus(\sigma_{2k-1})&=(A-\x{\sigma_{2k-1}},B+\x{\sigma_{2k-1}},C,D),\\
\sigma_{2k}&=(A-\x{\sigma_{2k-1}},B,C+\x{\sigma_{2k-1}},D).
\end{align*}
Now, $\sigma_{2k}\in\mdown_2$, so we have
$A-\x{\sigma_{2k-1}}=\x{\sigma_{2k}}$, or equivalently
$A=\{\x{\sigma_{2k-1}},\x{\sigma_{2k}}\}$.  Note, that
$\x{\sigma_{2k-1}}<\x{\sigma_{2k}}$.

By definition of $\mu$ this element is matched to
\[\muplus(\sigma_{2k})=(\x{\sigma_{2k}},B+\y{\sigma_{2k}},
C+\x{\sigma_{2k-1}}-\y{\sigma_{2k}},D).\]
Again, by our argument above, we  have
\[\sigma_{2k+1}=(\x{\sigma_{2k}}+x,B+\y{\sigma_{2k}}-x,
C+\x{\sigma_{2k-1}}-\y{\sigma_{2k}},D),\]
where $x$ is some element from $B+\y{\sigma_{2k}}$.

Note that $\x{\sigma_{2k}}<B$, by definition of the function $\x{}$,
and $\x{\sigma_{2k}}<\y{\sigma_{2k}}$, because
$\sigma_{2k}\in\mdown_2$.  It follows that $\x{\sigma_{2k}}<x$, so
$\x{\sigma_{2k}}=\x{\sigma_{2k+1}}$. Therefore $\sigma_{2k+1}$ is
matched to
\[\muplus(\sigma_{2k+1})=(x,
B+\y{\sigma_{2k}}-x+\x{\sigma_{2k}},
C+\x{\sigma_{2k-1}}-\y{\sigma_{2k}},D).\]

Finally, $\x{\sigma_{2k-1}}<\x{\sigma_{2k}}$ together with
$\x{\sigma_{2k}}=\x{\sigma_{2k+1}}$ implies that
$\x{\sigma_{2k-1}}<\x{\sigma_{2k+1}}$, for all $k$. This leads to
a~contradiction as we follow one turn of the cycle.  \qed

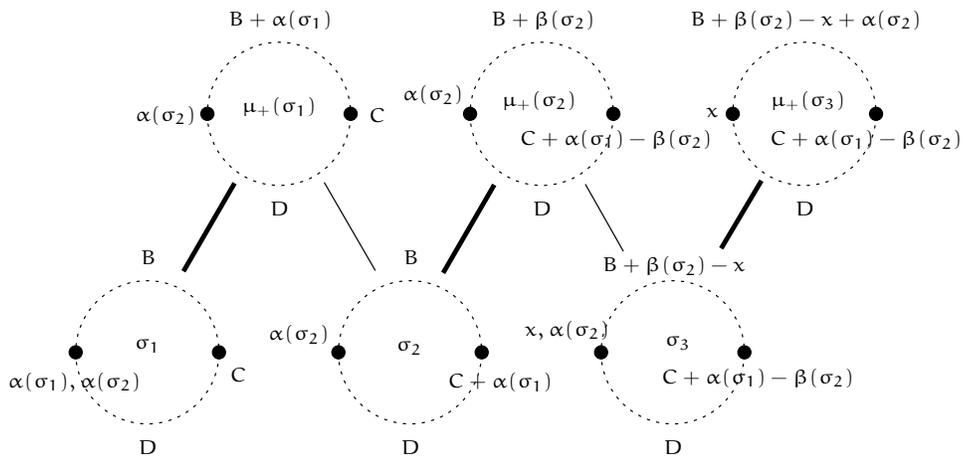
\myfig{acyc}{Closer look at the cycle at $\sigma_1$}


\section{Homology calculation}

In this section we shall apply the statements from \cref{sect:admt} to
our matching.  To start with, we need to extend our terminology.

\begin{df}
Let $\sigma$ be an arbitrary cube of $\Omega_n$. A {\bf
  $\Lambda$-path} starting at $\sigma$ is either $\sigma$ itself, if
$\sigma$ is critical, or a sequence of cubes
\begin{equation}
\label{eq:lp}
\sigma=\tau_1\prec\muplus(\tau_1)\succ\dots\succ\tau_q\prec\muplus(\tau_q)\succ\tau,
\end{equation}
where $q$ is a positive integer and $\tau$ is a critical
cube. Clearly, in the latter case, we must have $\sigma\in\mdown$.
\end{df}

We say that the $\Lambda$-path shown in \eqref{eq:lp} {\it ends at the cube
  $\tau$.} If the $\Lambda$-path consists of a single cube, we say
that this path {\it ends at $\sigma$.}

Note that removing the starting cube from an alternating path results
in a $\Lambda$-path. Likewise, when $q\geq 1$ in \eqref{eq:lp}, the
removal of the first two cubes from an alternating path, results in
a~new alternating path.

We shall next prove that the boundary operator in the Morse chain
complex is trivial, and the following lemmata will provide the key
building block of the argument.


\begin{lm} \label{lm:lp1}
Assume
$\sigma=\tau_1\prec\muplus(\tau_1)\succ\dots\succ\tau_q\prec\muplus(\tau_q)\succ\tau$
is a $\Lambda$-path, then $B(\tau_1)=\dots=B(\tau_q)=B(\tau)=\emptyset$.
\end{lm}
\pr 
By our construction and the description of the set of critical cubes,
we have the following facts:
\begin{enumerate}
\item[(1)] since $\tau$ is critical, and the dimension of $\tau$ is at most $n-3$,
we have $B(\tau)=\emptyset$;
\item[(2)] $|B(\tau_i)|+1=|B(\muplus(\tau_i))|$, for all $1\leq i\leq q$;
\item[(3)] the difference $|B(\muplus(\tau_i))|-|B(\tau_{i+1})|$ is
  either $0$ or $1$, for all $1\leq i\leq q$.
\end{enumerate}
It follows that, for all $1\leq i\leq q$, the difference
$|B(\tau_{i+1})|-|B(\tau_i)|$ equals to either $1$ or $0$, in
particular, $|B(\tau_i)|\leq|B(\tau_{i+1})|$.  Combining this
with the fact that $B(\tau)=\emptyset$, we can conclude that 
$B(\tau_i)=\emptyset$, for all $1\leq i\leq q$.
\qed

\begin{lm} \label{lm:lp2}
Assume $\sigma$ is a cube of $\Omega_n$, such that $\sigma\notin\mup$,
and $B(\sigma)=\emptyset$. Then, there exists a unique $\Lambda$-path
starting at $\sigma$. This path will end in a~critical cube
$\tau=(x,\emptyset,A(\sigma)+C(\sigma)-x,D(\sigma))$, where
$x=\max(A(\sigma)+C(\sigma))$.
\end{lm}
\pr If $\sigma$ is critical, then we must have $A(\sigma)=\x\sigma$,
and $\x\sigma>C(\sigma)$, so the statement of the lemma is correct.

Assume now that $\sigma\in\mdown$. By \cref{lm:lp1}, we know that in
any alternating path
\[\sigma=\tau_1\prec\muplus(\tau_1)\succ\dots\succ\tau_q\prec\muplus(\tau_q)\succ\tau\]
we must have $B(\tau_1)=\dots=B(\tau_q)=B(\tau)=\emptyset$.

We shall now prove the statement by  induction on $|A(\sigma)|$.
Consider first the basis case when $|A(\sigma)|=1$, say
$A(\sigma)=x$.  If $x>C(\sigma)$, then $\sigma$ is critical,
contradicting our assumption that $\sigma\in\mdown$. Therefore, we
must have $\y\sigma>x$. The only $\Lambda$-path starting from $\sigma$
is then
\begin{multline*}
\sigma\prec(x,\y\sigma,C(\sigma)-\y\sigma,D(\sigma))\succ
(x+\y\sigma,\emptyset,C(\sigma)-\y\sigma,D(\sigma))\\
\prec
(\y\sigma,x,C(\sigma)-\y\sigma,D(\sigma))\succ 
(\y\sigma,\emptyset,C(\sigma)-\y\sigma+x,D(\sigma)),
\end{multline*}
with the last cube being critical and satisfying the conditions of the lemma.
  
Next, we prove the induction step. Assume that
$A(\sigma)=\{x_1,\dots,x_k\}$, with $k\geq 2$, $x_1<\dots<x_k$, and
assume that the statement has been proved for smaller values of~$k$.
As we have shown in the proof of \cref{prop:clmacyc}, the only way a
$\Lambda$-path can start from $\sigma$ is
\[\sigma\prec(A(\sigma)-x_1,x_1,C(\sigma),D(\sigma))\succ
(A(\sigma)-x_1,\emptyset,C(\sigma)+x_1,D(\sigma)).\] We can now use
the induction assumption to conclude that there is a unique
$\Lambda$-path starting at the cube
$(A(\sigma)-x_1,\emptyset,C(\sigma)+x_1,D(\sigma))$.  This path will
end at the cube $(y,\emptyset,A(\sigma)+C(\sigma)-y,D(\sigma))$, where
$y=\max(A(\sigma)+C(\sigma))$.  Clearly, this is the unique
$\Lambda$-path starting at~$\sigma$.  \qed


We can now use the technical Lemmata \ref{lm:lp1} and \ref{lm:lp2} to
prove the main result of this section.

\begin{thm} \label{prop:bozero}
The boundary operator in the Morse chain complex is trivial.
\end{thm}
\pr Let $\sigma$ be a critical cube of dimension $d$, $d\geq 1$, say
$\sigma=(A,B,C,D)$.  Consider an alternating path
\[\sigma\succ\tau_1\prec\muplus(\tau_1)\succ\dots\succ\tau_q\prec\muplus(\tau_q)\succ\tau,\]
where $q$ is a nonnegative integer and $\tau$ is a critical cube of
dimension $d-1$.  For convenience, we set $\tau_{q+1}:=\tau$.  By
\cref{lm:lp1} we know that
$B(\tau_1)=\dots=B(\tau_q)=B(\tau)=\emptyset$.

Let us first consider the case when $B(\sigma)\neq\emptyset$.  In
particular, the cube $\sigma$ has dimension $n-2$. If $|B(\sigma)|\geq
2$, then $|B(\tau_1)|\geq 1$, contradicting the fact that $B(\tau_1)$
is empty. Therefore, in this case, there are no alternating paths
originating at $\sigma$ at all.

We can therefore assume that $B(\sigma)$ consists of a single element,
say $B(\sigma)=x$. We then have two options: either
$\tau_1=(A+x,\emptyset,C,D)$ or $\tau_1=(A,\emptyset,C+x,D)$.  By
\cref{lm:lp2} there exists a unique $\Lambda$-path starting from each
one. Both paths will end at the same critical cube
$(y,\emptyset,A+C+x-y,D)$, where $y=\max(A+C+x)$. By \cref{df:mc} the
boundary operator of the Morse complex evaluates to $0$ on~$\sigma$.

Let us now consider the case when $B(\sigma)=\emptyset$. Assume
$D(\sigma)=\{x_1,\dots,x_k\}$. Considering the alternating paths
starting from $\sigma$, there are $2k$ possibilities for the first
step.  The cube $\tau_1$ is either $(A+x_i,\emptyset,C,D-x_i)$ or
$(A,\emptyset,C+x_i,D-x_i)$, for $1\leq i\leq k$. By the argument
above, once $\tau_1$ is chosen, the rest of the alternating path can
be chosen uniquely. Accordingly we will have $2k$ alternating paths
starting at $\sigma$.  These paths will end at the cubes
$(y_i,\emptyset,A+C+x_i-y_i,D)$, where $y_i=\max(A+C+x_i)$, with
exactly two paths ending in each such cube. It follows from
\cref{df:mc} that the boundary operator of the Morse complex
evaluates to $0$ on this $\sigma$ as well.  \qed

We can now state the first main theorem of this paper.

\begin{thm} \label{crl:betti}
The Betti numbers of $\Omega_n$ are given by the following formulae:
\begin{align}
\beta_d&={n\choose d},\textrm{ for }0\leq d\leq n-3; \\
\beta_{n-2}&=2^n+{{n-1}\choose 2}-2.
\end{align}
\end{thm} 
\pr Combine \cref{thm:admt2}(1) with \cref{prop:count}.
\qed

The reader is invited to see how \cref{crl:betti} confirms the already
derived value for the M\"obius function of $\Omega_n$ via the
Euler-Poincar\'e formula .


\section{Explicit homology basis}


In this section we would like to describe an explicit homology basis
for $\Omega_n$.  To this end, let us define the following chains in
$C_*(\Omega_n)$.

\begin{df}
Assume $A$ and $C$ are two disjoint non-empty subsets of $[n]$.  Set
$S:=[n]\sm (A\cup C)$.  We define
\[
\rho_{A,C}:=\sum_{B\cup D=S} (A,B,C,D).
\]
\end{df}

Clearly, $\rho_{A,C}$ is a chain of dimension $|S|=n-|A|-|B|$.

\begin{prop} 
The chain $\rho_{A,C}$ is a cycle for all choices of $A$ and $C$.
\end{prop}
\pr 
Simple-minded boundary calculation gives
\begin{equation} 
\label{eq:rhob}
\begin{split}
\bo(\rho_{A,C})&=\bo(\sum_{B\cup D=S} (A,B,C,D))=
\sum_{B\cup D=S} \bo(A,B,C,D)\\
&=\sum_{B\cup D=S\atop x\in B}(A+x,B-x,C,D)+
   \sum_{B\cup D=S\atop x\in B}(A,B-x,C+x,D)\\
&+\sum_{B\cup D=S\atop x\in D}(A+x,B,C,D-x)+
   \sum_{B\cup D=S\atop x\in D}(A,B,C+x,D-x).
\end{split}
\end{equation}
All the terms on the right hand side of \cref{eq:rhob} are either of
the type $(A+x,B',C,D')$ or $(A,B',C+x,D')$, where $x\in S$, and $B'\cup
D'=S-x$.  Each term of the type $(A+x,B',C,D')$ occurs twice: once in
first sum, and once in the third one. Likewise, each term of the type $(A,B',C+x,D')$
occurs once in the second sum, and once in the third. In any case, the total
sum vanishes in $\dz$.
\qed

We can then see which critical cubes are contained in the support of this cycle.

\begin{prop}
\label{prop:rhocrit}
The presence of the critical cubes in the support of $\rho_{A,C}$ is
described by the following statements, where as before we set
$S:=[n]\sm (A\cup C)$.
\begin{enumerate}
\item[(1)] If $|A|\geq 2$, the support of $\rho_{A,C}$ contains no
critical cubes.
\item[(2)] If $A=a$ and $|C|\geq 2$, then the support of
  $\rho_{A,C}$ contains no critical cubes unless $a>C$. In the latter case,
it contains a unique critical cube $(a,\emptyset,C,S)$.
\item[(3)] If $A=a$, $C=c$, and $c-a\leq 1$, the support of
  $\rho_{A,C}$ contains a unique  critical cube $(a,\emptyset,c,S)$. 
\item[(4)] If $A=a$, $C=c$, and $c-a\geq 2$, the support of
  $\rho_{A,C}$ contains the critical cubes $(a,T,c,S\sm T)$, for any
  $T\subseteq S$, such that either $T=\es$ or $a<T<c$.
\end{enumerate}
\end{prop}
\pr This follows from applying case-by-case analysis to the
description of the critical cubes.  \qed

Note, that (2) and (3) of \cref{prop:rhocrit} means that we have the homology
generators dual to the cochains $(a,\emptyset,C,D)^*$, where 
\begin{enumerate}
\item[(1)] either $|C|\geq 2$, 
\item[(2)] or $C=c$ and $c-a\leq 1$.
\end{enumerate}
The first case covers all the homology in dimensions $n-3$ and lower.
The second case covers some of the homology in dimension $n-2$,
so that we are left to deal with the case $c-a\geq 2$.
This case is not quite as straightforward.


\begin{df}
Assume $G$ is an arbitrary subset of $[n]$, such that $|G|\geq 2$.
Set $a:=\min G$, $c:=\max G$.

Consider a decomposition of the set $[n]$ as a disjoint union
$E\coprod I\coprod G$, where $E=[n]\sm[a,c]$, and $I=[a,c]\sm G$. We
then set
\begin{equation*}
\gamma_G:=\sum_{i,j\in[a,c]}(i,B,j,D),
\end{equation*}
with two additional conditions on the summands: $D\subseteq I\cup E$ and $B\subseteq G\cup E$,
see \cref{fig:expl}.
\end{df}

Clearly, $\gamma_G$ is a chain for dimension $n-2$.

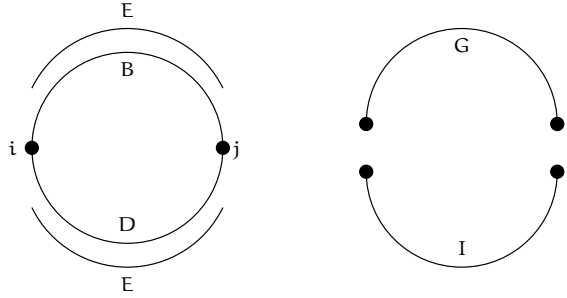
\myfig{expl}{Graphic description of the conditions on the labels of the cubes in the chain $\gamma_G$}

\begin{prop}
For every set $G$, such that $|G|\geq 2$, the chain $\gamma_G$ is a cycle.
\end{prop}
\pr
Just as in \cref{eq:rhob}, the straightforward boundary calculation yields
\begin{equation} 
\label{eq:gab}
\begin{split}
\bo(\gamma_G)&=\bo(\sum_{i,j\in[a,c]\atop {D\subseteq I\cup E\atop B\subseteq G\cup E}} (i,B,j,D))=
\sum_{i,j\in[a,c]\atop {D\subseteq I\cup E\atop B\subseteq G\cup E}} \bo(i,B,j,D)\\
&=\sum_{i,j\in[a,c]\atop {D\subseteq I\cup E\atop x\in B\subseteq G\cup E}}(i+x,B-x,j,D)+
   \sum_{i,j\in[a,c]\atop {D\subseteq I\cup E\atop x\in B\subseteq G\cup E}}(i,B-x,j+x,D)\\
&+\sum_{i,j\in[a,c]\atop {x\in D\subseteq I\cup E\atop B\subseteq G\cup E}}(i+x,B,j,D-x)+
   \sum_{i,j\in[a,c]\atop {x\in D\subseteq I\cup E\atop B\subseteq G\cup E}}(i,B,j+x,D-x).
\end{split}
\end{equation}
Let us analyze the right hand side of \cref{eq:gab}.  All the terms in
the sums are either of the type $(i+k,B',j,D')$ or $(i,B',j+k,D')$.
Due to symmetry, it is enough to consider the terms
$(i+k,B',j,D')$. These are subject to the conditions:
\begin{enumerate}
\item[(1)] $D'\subseteq I\cup E$ and $B'\subseteq G\cup E$,
\item[(2)] at least one of the elements $i$ and $k$ lies in $[a,c]=G\cup I$.
\end{enumerate}
Since we can swap $i$ and $k$, we can assume that either $i\in I$, or,
$i\in G$ and $k\notin I$. Label the sums on the right hand side of
\cref{eq:gab} with $I-IV$.  It is enough to show that each term
$(i+k,B',j,D')$ either appears once in two of these sums, or twice in
one of the sums.  The following table summarizes the 5 possible cases.
\begin{center}
\begin{tabular}{l | l} 
 $i$ and $j$ affiliations & sums where the term can be found   \\ [0.5ex] 
 \hline 
 $i,k\in I$ & twice in the sum III   \\  [0.5ex] 
 $i\in I$, $k\in G$ & in sums I and III  \\ [0.5ex] 
 $i\in I$, $k\in E$ & twice in the sum III \\ [0.5ex] 
 $i,k\in G$ & twice in the sum I   \\ [0.5ex] 
 $i\in G$, $k\in E$ & in sums I and III  \\ [0.5ex] 
\end{tabular}
\end{center}
\qed

\begin{prop}
The support of $\gamma_G$ contains a unique critical cube 
\[\sigma_G:=(\min G, G-\min G-\max G,\max G,I+E).\]
\end{prop}
\pr To start with, it is clear that $\sigma_G$ is critical and 
that it actually belongs to the support of $\gamma_G$.
Let $\tau=(i,B,j,D)$ be some other critical cube in the support
of $\gamma_G$. By definition of $\gamma_G$, we know that
$\min G,\max G\in B+i+j$. On the other hand, since $\tau$ is critical
we must have $i=\min(B+i+j)$ and $j=\max(B+i+j)$. It follows that
$i\leq\min G$ and $j\geq\max G$. On the other hand, we know that
$i,j\notin E$, so we must have  $i=\min G$ and $j=\max G$.

Finally, if $B\cap E\neq\emptyset$, then we cannot have both
$\min G=\min(B+i+j)$ and $\max G=\max(B+i+j)$. It follows, that $\sigma_G$
is the only critical cube in the support of~$\gamma_G$.
\qed

We finish by stating the second main theorem of this paper.

\begin{thm}
For $d=0,\dots,n-3$ a basis for the homology group $H_d(\Omega_n)$ is given by
the set $\myset{\rho_{a,C}}{a>C,\,n-|C|-1=d}$.

The basis for $H_{n-2}(\Omega_n)$ is given by the set 
\[\myset{\gamma_G}{|G|\geq 2}\cup\myset{\rho_{a,c}}{a>c}.\]
\end{thm}
\pr This follows from \cref{thm:admt2}(3).
\qed


\section{Final remarks}

It is not difficult to adopt our calculations to the case of other
coefficients. All one needs to do is to make a compatible choice of
incidence coefficients and then trace through our proofs, using the
general definition if the boundary operator in the Morse chain
complex.

Finally, we would like to remark that it would be interesting to
generalize our results to the case of more than $2$ points on
a~circle.



\begin{thebibliography}{AaA00}

\bibitem[BaK06]{BaK06} E.\ Babson, D.N.\ Kozlov, {\em Complexes of graph
homomorphisms},  Israel J.\ Math.\ {\bf 152} (2006), 285--312. 

\bibitem[BaK07]{BaK07} E.\ Babson, D.N.\ Kozlov, {\em Proof of the
    Lov\'asz Conjecture}, Annals of Math.\ (2) {\bf 165} (2007), no.\ 3, 
965--1007.

\bibitem[Ca09]{Ca} G.\ Carlsson, {\em Topology and data}, Bull.\ Amer.\ Math.\ 
Soc.\ (N.S.) {\bf 46} (2009), no.\ 2, 255--308.

\bibitem[EH10]{EH} H.\ Edelsbrunner, J.L.\ Harer, {\it Computational topology.
 An introduction.} American Mathematical Society, Providence, RI,
 2010.\ xii+241 pp.

\bibitem[Fu95]{Fu} W.\ Fulton, {\em Algebraic topology}, Graduate Texts in 
Mathematics {\bf 153}, Springer-Verlag, New York, 1995. xviii+430 pp.

\bibitem[GH81]{GH} M.J.\ Greenberg, J.R.\ Harper, {\em Algebraic Topology}, 
Mathematics Lecture Note Series {\bf 58}, Benjamin/Cummings Publishing
Co., Inc., Advanced Book Program, Reading, Mass., 1981. xi+311 pp.

\bibitem[Hat]{Hat} A.\ Hatcher, {\it Algebraic topology}, Cambridge
University Press, Cambridge, 2002.


\bibitem[Kn15]{Kn} K.P.\ Knudson, {\it Morse theory. Smooth and discrete.} 
World Scientific Publishing Co.\ Pte.\ Ltd., Hackensack, NJ, 2015.\
xiv+181 pp.

\bibitem[Ko07]{Ko07} D.N.\ Kozlov, 
{\em Chromatic numbers, morphism complexes, and Stiefel-Whitney
characteristic classes}, in: {\em Geometric Combinatorics} (eds.\ E.\
Miller, V.\ Reiner, B.\ Sturmfels), pp.\ 249--315, IAS/Park City
Mathematics Series {\bf 13}, American Mathematical Society,
Providence, RI; Institute for Advanced Study (IAS), Princeton, NJ;
2007.

\bibitem[Ko05]{admt} D.N.\ Kozlov, \textit{Discrete Morse theory for free chain complexes.} 
(English, French summary) C.\ R.\ Math.\ Acad.\ Sci.\ Paris {\bf 340} (2005), no.\ 12, 867--872.

\bibitem[Ko08]{Ko08} D.N.\ Kozlov, \textit{Combinatorial Algebraic Topology},
Algorithms and Computation in Mathematics {\bf 21}, Springer, Berlin,
2008, xx+389 pp.

\bibitem[Ko20a]{Ko20a} D.N.\ Kozlov, \textit{A combinatorial method to
  compute explicit homology cycles using discrete Morse
  theory},  J.\ Appl.\ Comput.\ Topol.\ {\bf 4} (2020), no.\ 1,
  79--100.

\bibitem[Ko20b]{Ko20b} D.N.\ Kozlov, \textit{Organized collapse: an
  introduction to discrete Morse theory}, Graduate Studies in
  Mathematics {\bf 207}, American Mathematical Society, Providence,
  RI, 2020, xxiii+312 pp.

\bibitem[Ko21]{Ko21} D.N.\ Kozlov, {\it Stirling complexes}, 2021, preprint 13 pp.,
to appear in Journal of Applied and Computational Topology.

\bibitem[MKS]{MKS} W.\ Magnus, A.\ Karrass, D.\ Solitar, {\it
  Combinatorial group theory.  Presentations of groups in terms of
  generators and relations.} Reprint of the 1976 second edition. Dover
  Publications, Inc., Mineola, NY, 2004. xii+444 pp.

\bibitem[Mu84]{Mu} J.R.\ Munkres, {\em Elements of algebraic topology},
Addison-Wesley Publishing Company, Menlo Park, CA, 1984 ix+454 pp.

\bibitem[Sc19]{Sc19} 
N.\ Scoville, {\em Discrete Morse theory}, Student Mathematical
Library 90, American Mathematical Society, 273 pp.

\end{thebibliography}
\end{document}